\magnification=1200

\def\C{{\bf C}}
\def\F{{\cal F}}
\def\hal{{\vrule height 10pt width 4pt depth 0pt}}
\def\K{{\cal K}}
\def\la{{\langle}}

\def\PK{{\cal PK}}
\def\R{{\bf R}}
\def\ra{{\rangle}}
\def\Z{{\bf Z}}
\centerline{\bf Noncommutative complex analysis and}

\centerline{\bf Bargmann-Segal multipliers}
\bigskip

\centerline{Richard Rochberg and Nik Weaver}
\bigskip
\bigskip

{\narrower{
\noindent \it We state several equivalent noncommutative versions of the
Cauchy-Riemann equations and characterize the unbounded operators on
$L^2(\R)$ which satisfy them. These
operators arise from the creation operator via a functional calculus
involving a class of entire functions, identified by Newman and Shapiro
([13], [14]), which act as unbounded multiplication operators on
Bargmann-Segal space.
\bigskip}}
\bigskip

\noindent {\bf 1. Noncommutative Cauchy-Riemann equations.}
\bigskip

Various branches of the field of operator algebras are, following
Connes [6], Effros [7], Woronowicz [17], and others, often
viewed as ``noncommutative'' or ``quantized'' versions of classical
mathematical subjects. Underlying this attitude is a rather diverse family of
fundamental examples, of which the noncommutative tori [15] and the
noncommutative plane ([11], [12]) are probably the most basic.
Our aim here is to determine an appropriate version of complex structure
on the latter.
\medskip

The basic elements of the noncommutative plane are these. The underlying
Hilbert space is $L^2(\R)$. Compact operators play the role of continuous
functions which vanish at infinity and bounded operators play the role of
bounded measurable functions. The unbounded self-adjoint operators $Q = M_x$
(multiplication by $x$) and $P = -i\hbar{d\over{dx}}$, where $\hbar$ is a
nonzero real number, are analogous to the coordinate functions, and commutation
with the operators $-{1\over\hbar}P$ and ${1\over\hbar}Q$ correspond to the
partial derivatives ${\partial\over{\partial x}}$ and
${\partial\over{\partial y}}$.
\medskip

Physically, this structure characterizes a spinless, nonrelativistic,
one-dimensional quantum-mechanical particle. The ``noncommutative plane''
intuition is motivated by the fact that
the phase space of a classical one-dimensional particle is the ordinary plane,
the classical position and momentum observables are the coordinate functions
on $\R^2$, and so on.
\medskip

In what sense does the noncommutative plane carry a complex structure?
Consider a differentiable function $\phi: \C \to \C$. By regarding it as
a function from $\R^2$ into $\C$, we may condense the Cauchy-Riemann
equations, which diagnose whether $\phi$ is holomorphic, into the single
condition
$${{\partial \phi}\over{\partial y}}  = i{{\partial \phi}\over{\partial x}}.$$
In terms of the unbounded multiplication operator $M_\phi: f \to \phi f$
acting on $L^2(\R^2)$, this condition can formally be reexpressed as
$$[M_\phi, {\partial\over{\partial x}} + i{\partial\over{\partial y}}] = 0$$
or equivalently
$$[M_\phi, {\partial\over{\partial y}} - i{\partial\over{\partial x}}] = 0.$$
(The general principle is that $[M_\phi, {\partial\over{\partial x}}]
= -M_{\partial \phi/\partial x}$, etc.)
\medskip

Now the noncommutative analogs of ${\partial\over{\partial x}}$ and
${\partial\over{\partial y}}$ are commutation with $-{1\over\hbar}P$ and
${1\over\hbar}Q$. So, depending on whether we
identify position with the $x$-axis and momentum with the $y$-axis or
vice versa, the noncommutative version of the Cauchy-Riemann
equations is either the condition
$$[A, a^-] = 0\eqno{(*)}$$
or the condition
$$[A, a^+] = 0,\eqno{(**)}$$
where $a^-$ is the annihilation operator $a^- = {1\over 2}(Q + iP)$ and
$a^+$ is the creation operator $a^+ = {1\over 2}(Q - iP)$. Either condition
is to be satisfied by a possibly unbounded operator $A$ acting on $L^2(\R)$.
\medskip

In classical complex analysis the variables can be switched without
essentially changing the theory, but in the noncommutative setting this
is not true: the creation and annihilation operators are
not equivalent. However, if one uses the right domains they
are each other's adjoints, so the
class of operators which satisfy ($*$) ought to be adjoint to the class of
operators which satisfy ($**$).
\medskip

This is not a rigorous assertion since ($*$) and ($**$) are only formal
expressions; there is no standard way to define the commutator of two
unbounded operators. To make these conditions precise we need a clearer
picture of the operators $a^+$ and $a^-$.
Mathematically, the nicest representation of the
creation and annihilation operators --- and the setting in which their natural
domains are apparent --- is on Bargmann-Segal space, so we go there to
pursue the question further.
\bigskip
\bigskip

\noindent {\bf 2. Bargmann-Segal formulation.}
\bigskip

From now on we set $r = 1/\hbar$ and assume $r > 0$. (Changing the sign of
$\hbar$ is tantamount to interchanging $P$ and $Q$.)
The Bargmann-Segal space, also known as Bargmann-Fock or Fischer space,
is the Hilbert space $\F$ consisting of all entire
functions on the complex plane which are square-integrable with respect
to the Gaussian measure $\mu$ defined by $d\mu(z) =
(r/\pi) e^{-r|z|^2}dz$. It is isomorphic to $L^2(\R)$ via the
Bargmann-Segal transform, which takes the natural basis of $\F$ --- the
orthogonal functions $z^n$, appropriately normalized --- to the Hermite
basis of $L^2(\R)$. Under this transformation the operators $a^+$ and $a^-$
on $L^2(\R)$ are respectively identified with the operators $M_z$ and
${1\over r}{d\over{dz}}$ on $\F$. The natural domain of $M_z$ is
$\{f \in \F: zf \in \F\}$ and the corresponding domain of
${1\over r}{d\over{dz}} = M_z^*$ is $\{f \in \F: \bar z f \in L^2(\C, d\mu)\}$
(see e.g.\ [5]). See [8] for a thorough treatment of the Bargmann-Segal
transform and see [2], [9], [10] for work on Toeplitz operators on $\F$.
\medskip

Let $\K$ be the linear span in $\F$ of the reproducing kernels
$e_w(z) = e^{rz\bar w}$ ($w \in \C$). These functions have the property
that $\la f, e_w\ra = f(w)$ for all $f \in \F$.
Let $\PK$ be the algebra generated
by $\K$ together with the function $z$. Also, let $\Lambda$ be the
``operator class'' of Newman and Shapiro ([13], [14]); this consists of
those entire functions $\phi$ which satisfy $|\phi(z)|
= O(e^{r|z|^2/2 - N|z|})$ for
all $N > 0$. The following result is given in the unpublished manuscript
[14]; for the reader's convenience we include the proof.
\bigskip

\noindent {\bf Theorem 1.} ([14], $\S$ 2.2.) {\it Let $\phi$ be an
entire function on $\C$. Then $\phi \in \Lambda$ if and only if
$\phi e_w \in \F$ for all $w \in \C$. The class $\Lambda$ is closed
under differentiation.}
\medskip

\noindent {\it Proof.} The forward direction of the first assertion
is a routine calculation.
For the reverse direction, suppose $\phi e_w \in \F$ for all $w$
and fix $N > 0$. Let $w_k = (\sqrt{2}N/r)e^{2\pi i(k/4)}$ for
$0 \leq k \leq 3$. Then
$$|\phi(z)e^{r z\bar w_k}| = |\la \phi e_{w_k}, e_z\ra|
\leq \|\phi e_{w_k}\|e^{r|z|^2/2}$$
for each $k$, by the Cauchy-Schwartz inequality. Also, for any $z \in \C$
there is a value of $k$ such that $-\pi/4 \leq {\rm Arg}(z\bar w_k)
\leq \pi/4$, and for this $k$ we have
${\rm Re}(z\bar w_k) \geq |z\bar w_k|/\sqrt{2} = N|z|/r$. Thus
$$|\phi(z)e^{r z\bar w_k}| = |\phi(z)|e^{r {\rm Re}(z\bar w_k)}
\geq |\phi(z)|e^{N|z|}$$
and we conclude that
$$|\phi(z)| \leq C\cdot e^{r|z|^2/2 - N|z|}$$
where $C = {\rm max}_{0 \leq k \leq 3} \|\phi e_{w_k}\|$. This shows
that $\phi \in \Lambda$.
\medskip

Closure of $\Lambda$ under differentiation follows from the formula
$$\eqalign{\phi'(w) &= \la {{d}\over{d z}}\phi, e_w\ra
= \la \phi, rze_w\ra\cr
&= (r^2/\pi)\int \bar{z}\phi(z)e^{r(\bar{z}w - |z|^2)}dz;\cr}$$
the bound $|\phi(z)| = O(e^{r|z|^2/2 - (N + 1)|z|})$ then yields
$$\eqalign{|\phi'(w)|
&\leq C\int |z|e^{r|z|^2/2 - (N+1)|z|}\,
e^{r({\rm Re}(z\bar{w}) - |z|^2)}\, dz\cr
&= C\int |z|e^{r|w|^2/2}e^{-r|z-w|^2/2}\, e^{-(N+1)|z|}\,dz\cr
&= Ce^{r|w|^2/2 - N|w|}
\int |z|e^{-r|z-w|^2/2 - N(|z| - |w|)}\, e^{-|z|}\, dz.\cr}$$
The final integral is finite because $|z-w| \geq \big{|}|z| - |w|\big{|}$,
so that the first exponent in the integrand is bounded above.
Thus $\phi'(w) = O(e^{r|w|^2/2 - N|w|})$, as desired.\hfill\hal
\bigskip

We remark that $\Lambda$ is also closed under integration; this follows
from an elementary estimate based on the inequality
$$|\phi(z)| \leq |\phi(0)| + \int_0^{|z|} |\phi'(t{z\over{|z|}})|dt.$$

For $\phi \in \Lambda$, let $M_\phi$ be the operator of multiplication by
$\phi$, with domain $D(M_\phi) = \{g \in \F: \phi g \in \F\}$.
Thus $\K \subset D(M_\phi)$
by the preceding result, and in fact $\PK \subset D(M_\phi)$ since every
function in $\PK$ is of exponential type, i.e.\ is $O(e^{N|z|})$ for some
$N$. The adjoint $M_\phi^*$ has similar properties:
\bigskip

\noindent {\bf Proposition 2.} {\it Let $\phi \in \Lambda$. Then
$\PK \subset D(M_\phi^*)$, $\K$ is a core for $M_\phi^*$, and we have
$M_\phi^*(e_w) = \bar{\phi}(w)e_w$ for all $w \in \C$.}
\medskip

\noindent {\it Proof.} For $w \in \C$ and $g \in D(M_\phi)$ we have
$$\la M_\phi g, e_w\ra = \phi(w)g(w) = \la g, \bar{\phi}(w)e_w\ra .$$
This shows that $\K \subset D(M_\phi^*)$ and $M_\phi^*(e_w) =
\bar{\phi}(w)e_w$.
\medskip

Let $A$ be the closure of $M_\phi^*|_\K$. To prove that $\K$ is a core for
$M_\phi^*$, we must show that $M_\phi^* \subset A$, or equivalently that
$A^* \subset M_\phi$. Let $g \in D(A^*)$ and let $h = A^*g$.
Then we have
$$h(w) = \la h, e_w\ra = \la g, A e_w\ra = \la g, \bar{\phi}(w)e_w\ra
= \phi(w)g(w)$$
for all $w \in \C$. Thus $\phi g = h \in \F$, hence $g \in D(M_\phi)$,
and so $A^* \subset M_\phi$ as desired.
\medskip

Finally, $\PK \subset D(M_\phi^*)$ follows from the standard fact that
$D(M_\phi^*)$ contains $\{f \in \F: \bar\phi f \in L^2(\C, d\mu)\}$.
If $f \in \PK$ then $\bar\phi f$ satisfies the same
growth condition as $\phi$ and hence is square-integrable, so this
implies $f \in D(M_\phi^*)$.\hfill\hal
\bigskip

We will argue in the next section that the operators $M_\phi$ for
$\phi \in \Lambda$ are the natural solutions to ($**$) and their adjoints
are the natural solutions to ($*$).
\bigskip
\bigskip

\noindent {\bf 3. The main result.}
\bigskip

The following seems to be a reasonable notion of commutation for unbounded
operators. In our setting the classes $\K$ and $\PK$ will play the role of $E$.
\bigskip

\noindent {\bf Definition 3.} Let $E$ be a linear subspace of a Hilbert space
$H$ and let $A$ and $B$ be unbounded operators on $H$ both of whose domains
contain $E$, that is $E \subset D(A) \cap D(B)$. Suppose also that
$E \subset D(A^*) \cap D(B^*)$. Then we say that $A$ and $B$
{\it commute relative to $E$} and write $[A, B] = 0$ (rel $E$) if
$$\la A\eta,B^*\xi\ra = \la B\eta, A^*\xi\ra$$
for all $\eta, \xi \in E$.\hfill\hal
\bigskip

We now present our main theorem. We regard conditions (a) and (b) as justifying
the idea that operators which satisfy these conditions are the natural
solutions to equation ($**$): (b) is the stronger version, which says that
$A$ commutes with $M_z$ relative to the ``large'' domain $\PK$, while
(a) involves the ``small'' domain $\K$ and the morally weaker condition that
$A$ commutes with $M_{e_w}$ for all $w \in \C$. Likewise, conditions (f)
and (g) are weak and strong versions of equation ($*$).
\bigskip

\noindent {\bf Theorem 4.} {\it Let $A$ be a closed, densely defined,
unbounded operator on $\F$. Suppose $\K \subset D(A) \cap D(A^*)$ and
$K$ is a core for $A^*$. Then the following are equivalent:
\medskip

{\narrower{
\noindent (a) $[A, M_{e_w}] = 0$ (rel $\K$) for all $w \in \C$;
\smallskip

\noindent (b) $\PK \subset D(A) \cap D(A^*)$ and $[A, M_z] = 0$ (rel $\PK$);
\smallskip

\noindent (c) $A = M_\phi$ for some $\phi \in \Lambda$;
\smallskip

\noindent (d) $A^*e_w = \bar{\phi}(w)e_w$ for some $\phi \in \Lambda$
and all $w \in \C$;
\smallskip

\noindent (e) $e_w$ is an eigenvector of $A^*$ for all $w \in \C$;
\smallskip

\noindent (f) $[{d\over{dz}}, A^*] = 0$ (rel $\K$); and
\smallskip

\noindent (g) $\PK \subset D(A) \cap D(A^*)$ and
$[{d\over{dz}}, A^*] = 0$ (rel $\PK$).
\medskip}}}

\noindent {\it Proof.}
\medskip

(a) $\Rightarrow$ (c). Suppose (a) holds and let
$\phi = A(e_0) = A(1)$. Then for any $v,w \in \C$ we have
$$\la Ae_v, e_w\ra = \la M_{e_v}e_0, A^*e_w\ra
= \la Ae_0, M_{e_v}^*e_w\ra = e_v(w)\la \phi, e_w\ra = \phi(w)e_v(w);$$
this implies that $Ae_v = \phi e_v$. Hence $\phi e_v \in \F$ for all $v$,
so $\phi \in \Lambda$ follows from Theorem 1. Also $A|_\K = M_\phi|_\K$.
\medskip

Let $A_0 = A|_\K$. Then the above shows that $A_0 \subset M_\phi$ and
therefore $M_\phi^* \subset A_0^*$. But also $A^* \subset A_0^*$, and $\K$
is a core for both $A^*$ (by hypothesis) and $M_\phi^*$ (by Proposition 2). So
$A^* = M_\phi^*$ and hence $A = M_\phi$.
\medskip

(c) $\Rightarrow$ (d). This was shown in Proposition 2.
\medskip

(d) $\Rightarrow$ (e). Vacuous.
\medskip

(e) $\Rightarrow$ (a). Suppose every $e_w$ is an eigenvector of $A^*$ and
fix $u,v,w \in \C$. Let $A^*e_w = \bar{\lambda} e_w$. Then
$$\la M_{e_u}e_v, A^*e_w\ra
= \lambda \la e_{u + v}, e_w\ra
= \lambda e^{r(\bar u + \bar v)w}$$
while
$$\la Ae_v, M_{e_u}^* e_w\ra
= \la Ae_v, e^{ru\bar w}e_w\ra
= e^{r\bar u w}\la e_v, A^*e_w\ra
= \lambda e^{r\bar u w}\la e_v, e_w\ra
= \lambda e^{r\bar u w}e^{r\bar v w}.$$
So $[A, M_{e_u}] = 0$ (rel $\K$).
\medskip

(c) $\Rightarrow$ (b). $\PK \subset D(M_\phi)$ is easy and $\PK \subset
D(M_\phi^*)$ was shown in Proposition 2. If $\eta, \xi \in \PK$ then
$$\la M_\phi \eta, M_z^*\xi\ra = \la \phi\eta, M_z^*\xi\ra
= \la z\phi\eta, \xi\ra$$
and
$$\la M_z\eta, M_\phi^*\xi\ra = \la z\eta, M_\phi^*\xi\ra
= \la \phi z\eta, \xi\ra.$$
So $[M_\phi, M_z] = 0$ (rel $\PK$).
\medskip

(b) $\Rightarrow$ (g). This is immediate from the fact that
$M_z^* = {1\over r}{d\over{dz}}$ and the definition of commutation.
\medskip

(g) $\Rightarrow$ (f). Vacuous.
\medskip

(f) $\Rightarrow$ (e). Suppose $[{d\over{dz}}, A^*] = 0$ (rel $\K$).
Then for any $v,w \in \C$
$$\la A^*e_v, ze_w\ra
= \la {1\over r}{d\over{dz}}e_v, Ae_w\ra
= \la \bar v e_v, Ae_w\ra
= \la A^*e_v, ve_w\ra.$$
It follows that $A^*e_v$ is orthogonal to $(z-v)f$
if $f$ is any finite linear combination of the functions $e_w$.
\medskip

In particular let
$$f_h = {1\over{h^n}} \sum_{k=1}^n (-1)^{n-k} {n\choose k} e_{kh}.$$
Then $f_h$ is an $n$th order difference quotient and it converges in $\F$ to
the $n$th Hermite polynomial $H_n$ as $h \to 0$. We also
have $(z-v)f_h \to (z-v)H_n$ in $\F$. As the $H_n$ span all polynomials,
it follows that $Ae_v$ is orthogonal to every polynomial which
vanishes at $v$. So letting $\lambda = \la Ae_v, e_v\ra/\|e_v\|^2$, we have
that $Ae_v - \lambda e_v$ is orthogonal to every polynomial. Hence
$Ae_v = \lambda e_v$, as desired.\hfill\hal
\bigskip

Thus, if we fix the convention that ($**$) is the noncommutative version of
the Cauchy-Riemann equations, then
the unbounded operators which play the role of ``noncommutative''
entire functions are the operators $M_\phi$ ($\phi \in \Lambda$) --- which
is somewhat surprising considering that these operators all commute.
Noncommutativity resides in the fact that they are not normal.
\medskip

In any reasonable sense $M_\phi = \phi(M_z)$, so
these are just the operators that one obtains via functional calculus by
applying functions in the operator class $\Lambda$ to the creation
operator. Thus this class of operators possesses natural
linear, algebraic, and topological structures isomorphic to those on
$\Lambda$.
\bigskip

\noindent {\bf Remark 5.}
It is also possible to differentiate these operators with respect to the
two noncommutative coordinates, in the sense of taking the commutator
with $-{1\over \hbar}P$ and ${1\over\hbar}Q$.
Since $Q = (a^+ + a^-)$ and $P = i(a^+ - a^-)$, in
the Bargmann-Segal representation these operators are given by
$Q = (M_z + {1\over r}{d\over{dz}})$ and
$P = i(M_z - {1\over r}{d\over{dz}})$. So if $\phi \in \Lambda$
and $\eta \in \F$ is of exponential type then
$$\eqalign{[{1\over\hbar}Q, M_\phi]\eta
&= {1\over\hbar}Q(M_\phi\eta) - M_\phi({1\over\hbar}Q\eta)\cr
&= r(Q(\phi\eta) - \phi\cdot Q\eta)\cr
&= r( z\phi\eta + {1\over r}(\phi'\eta + \phi\eta')
- z\phi\eta -{1\over r}\phi\eta')\cr
&= M_{\phi'}\eta\cr}$$
and
$$\eqalign{[-{1\over\hbar}P, M_\phi]\eta
&= -{1\over\hbar}P(M_\phi\eta) - M_\phi(-{1\over\hbar}P\eta)\cr
&= -r(P(\phi\eta) - \phi\cdot P \eta)\cr
&= -ir(z\phi\eta - {1\over r}(\phi'\eta + \phi\eta')
- z\phi\eta +{1\over r}\phi\eta')\cr
&= M_{i\phi'}\eta.\cr}$$
Thus differentiation of $M_\phi$ in the noncommutative sense corresponds
to differentiation of the symbol $\phi$, and we have also verified that
$M_\phi$ satisfies the Cauchy-Riemann condition in the form
$$[-{1\over\hbar}P, M_\phi] = i[{1\over\hbar}Q, M_\phi].$$

\noindent {\bf Remark 6.} We can also identify a reasonable noncommutative
version of harmonic functions. Here the classical condition is
${{d^2 \phi}\over{dzd\bar{z}}} = 0$, which we can
formally quantize as
$$[a^+, [a^-, A]] = [M_z, [M_z^*, A]] = 0.\eqno{(\dag)}$$
The simplest way to make this harmonicity condition precise is to invoke
Theorem 4 to interpret it as asserting that $[M_z^*, A] = M_\phi$ for
some $\phi \in \Lambda$, and then to require, minimally, that
$$M_z^*Ae_v - AM_z^*e_v = M_\phi e_v$$
for all $v \in \C$. (In particular, we ask that $Ae_v \in D(M_z^*)$ for all
$v \in \C$; a short computation
shows that this is equivalent to the requirement that $e_v \in D([M_z^*, A])$
in the sense that the map $\eta \mapsto \la Ae_v, M_z\eta\ra -
\la M_z^*e_v, A^*\eta\ra$ is bounded on $\K$.)
\medskip

On $\K$, the solutions to ($\dag$) are precisely those of the
form $A = M_\Phi + B$ where $\Phi' = \phi$ and $B$ satisfies
$[M_z^*, B] = 0$ (rel $\K$). So by Theorem 4 (f), $B|_\K = M_\Psi^*|_\K$
for some $\Psi \in \Lambda$. Thus, ($\dag$) implies that
$A = M_\Phi + M_\Psi^*$ on $\K$, for some $\Phi, \Psi \in \Lambda$.
\medskip

In fact, if we assume $\K \subset D(A) \cap D(A^*)$ then $\K$ must be a core
for $A$. For letting $A_1$ be the closure of $(M_\Phi + M_\Psi^*)|_\K$ and
$A_2$ the closure of $(M_\Phi^* + M_\Psi)|_\K$, we have $A_2 \subset A_1^*$
and $A_1 \subset A_2^*$, hence $A_1^* = A_2$; so $A_2 \subset A^*$ implies
$A \subset A_1$. We conclude that reasonable domain restrictions
imply that ``noncommutative harmonic operators'' have $\K$ as a core and are
of the form $M_\Phi + M_\Psi^*$ ($\Phi, \Psi \in \Lambda$)
on $\K$. Conversely, a short computation
shows that any operator of this form does satisfy ($\dag$).
\medskip

If we instead begin with $[M_z^*, [M_z, A]] = 0$ as the defining condition
of harmonicity,
taking adjoints yields $[M_z, [M_z^*, A^*]] = 0$, so that (again, up to
reasonable domain
assumptions) $A^*$ is of the form $M_\Phi + M_\Psi^*$, and hence so is $A$.
\bigskip
\bigskip

\noindent {\bf 4. Counterexamples.}
\bigskip

The requirement that $\K$ be a core for $A^*$ in Theorem 4 is probably
unnecessary. It was conjectured in [13] (see also [4], [5])
that $\K$ is a core for every
$M_\phi$ ($\phi \in \Lambda$), and if this is true then our hypothesis that
$\K$ is a core for $A^*$ is superfluous. For the only place that it is
used is in the proof of (a) $\Rightarrow$ (c), and if we knew that $\K$
was a core for $M_\phi$ we could reason instead as follows: we know
$A|_\K = M_\phi|_\K$, hence $M_\phi \subset A$ since $\K$ is a core for
$M_\phi$; therefore $A^* \subset M_\phi^*$, but $\K \subset D(A^*)$ and
$\K$ is also a core for $M_\phi^*$, hence $A^* = M_\phi^*$ and so
$A = M_\phi$.
\medskip

On the other hand, the requirement $\K \subset D(A) \cap D(A^*)$ is crucial.
There are operators for which this is not true but which do commute with
$M_z$ is a reasonable sense, namely any multiplication operator $M_\phi$
with $\phi \not\in \Lambda$.
\medskip

Moving to this broader setting introduces a number of pathologies.
First, observe that since $M_z$ is unbounded there must exist $f \in \F$
such that $zf \not \in \F$. Then $1 \in D(M_f)$ but $z \not\in D(M_f)$.
In fact, using the fact that the norm of $z^n$ in $\F$ is $r^{n/2}\sqrt{n!}$
it is easy to check that the function $f(z) = \sum a_n z^n$ with
$a_n = (r^n n! (n+1)^2)^{-1/2}$ has this property, for example. Furthermore,
the function $g(z) = \sum a_{n+k} z^n$ satisfies $pg \in \F$ if $p$ is
a polynomial of degree at most $k$, and $pg \not\in \F$ if $p$ is a
polynomial of degree larger than $k$.
\medskip

Note that in general if $D(M_\phi)$ contains any polynomial $p$ of degree
$k$ then it contains all polynomials of degree $\leq k$. This follows
immediately from the definitions and the fact that at $\infty$,
$|p| \sim c|r|^{2n}$.
\medskip

In fact, given any nonnegative integer $k$ it is possible to find an
entire $\phi$ so that $D(M_\phi)$ consists exactly of ${\cal P}_k =
\{f: f$ is a polynomial of degree $\leq k\}$. A convenient way to see
this is to follow the lead of Seip and Wallstein. In Proposition 2.1 of
[16] they show, in our language, that if $\phi$ is the Weierstrass
$\sigma$-function with period lattice $\{\sqrt{\pi\over r}(n + im):
n, m \in \Z\}$ then $D(M_\phi) = \{0\}$. Suppose now that $k$ is given.
Select any $k + 2$ of the zeros of $\phi$ and let $p$ be a polynomial
with exactly those points as its (simple) zeros. The natural modification
of the argument in [16] gives
$$D(M_{\phi/p}) = \Big\{f: f\hbox{ entire and }
\int_{|z| > 1} |f|^2 |z|^{-2(k+2)} dz < \infty\Big\}.$$
This last set is ${\cal P}_k$, as required. (We note, as is pointed out in
[16], that although this particular choice of $\phi$ is extremely convenient
computationally it is not essential for such an argument.)
\medskip

It seems reasonable to restrict attention to the case that $M_\phi$ has
dense domain, but there are still serious problems in this case. For any
functions $f$ and $g$ of exponential type, any $w \in \C$ with $|w| < 1$,
and any $a \in (1-{1\over{2|w|}}, {1\over{2|w|}})$
we have $g e^{-awz^2} \in \F$ and
$$(f e^{wz^2})(g e^{-awz^2}) = fg e^{(1-a)wz^2} \in \F$$
since $|-aw| < {1\over 2}$ and $|(1-a)w| < {1\over 2}$.
It follows that $M_\phi$, where $\phi = f e^{wz^2}$, has a ``large'' domain,
but if $|w| \geq {1\over 2}$ then
nothing in its domain has exponential type, and in particular $\PK$ is disjoint
from $D(M_\phi)$. From a comment on
page 126 of [3] it appears that there are similar examples of such $\phi$ of
order 2 and arbitrarily large type, such that $D(M_\phi)$ contains functions
of order 2 and type
less than ${1\over 2}$. Another interesting example is given in [4]; this
is a function $\phi$ such that every polynomial belongs to $D(M_\phi)$ but
$\K \not\subset D(M_\phi)$ and the polynomials are not a core of $M_\phi$.
\medskip

For these reasons, from a purely operational point of view it clearly is
desirable to require that $\K \subset D(A) \cap D(A^*)$. Some
philosophical motivation may also come from the view that the coherent
states $e_w$ play the role of the points of the noncommutative plane
(see e.g.\ [1]).
\bigskip
\bigskip

We wish to thank Al Baernstein,
James Jenkins, John McCarthy, and Harold Shapiro for helpful discussions.
\bigskip
\bigskip

[1] S.\ T.\ Ali, J.-P.\ Antoine, J.-P.\ Gazeau, and U.\ A.\ Mueller,
Coherent states and their generalizations: a mathematical overview,
{\it Rev.\ Math.\ Phys.\ \bf 7} (1995), 1013-1104.
\medskip

[2] C.\ A.\ Berger and L.\ A.\ Coburn, Toeplitz operators on the Segal-Bargmann
space, {\it Trans.\ Amer.\ Math.\ Soc.\ \bf 301} (1987), 813-829.
\medskip

[3] R.\ Boas, {\it Entire Functions}, Academic Press (1954).
\medskip

[4] A.\ A.\ Borichev, The polynomial approximation property in Fock-type
spaces, {\it Math.\ Scand.\ \bf 82} (1998), 256-264.
\medskip

[5] D.\ Cicho\'n and J.\ Stochel, On Toeplitz operators in Segal-Bargmann
spaces, {\it Univ.\ Iagel.\ Acta Math.\ \bf 34} (1997), 35-43.
\medskip

[6] A.\ Connes, {\it Noncommutative geometry}, Academic Press (1995).
\medskip

[7] E.\ G.\ Effros, Advances in quantized functional analysis,
{\it Proc.\ ICM (Berkeley, California, 1986)}, 906-916, AMS (1987).
\medskip

[8] G.\ B.\ Folland, {\it Harmonic Analysis in Phase Space},
Princeton (1989).
\medskip

[9] J.\ Janas, Unbounded Toeplitz operators in the Bargmann-Segal space,
{\it Studia Math.\ \bf 99} (1991), 87-99.
\medskip

[10] J.\ Janas and J.\ Stochel, Unbounded Toeplitz operators in the
Segal-Bargmann space, II, {\it J.\ Funct.\ Anal.\ \bf 126} (1994), 418-447.
\medskip

[11] J.\ E.\ Moyal, Quantum mechanics as a statistical theory,
{\it Proc.\ Cambridge Philos.\ Soc.\ \bf 45} (1949), 99-124.
\medskip

[12] J.\ von Neumann, Die eindeutigkeit der Schr\"odingerschen
operatoren, {\it Math.\ Ann.\ \bf 104} (1931), 570-578.
\medskip

[13] D.\ J.\ Newman and H.\ S.\ Shapiro, Fischer spaces of entire functions,
in {\it Entire Functions and Related Parts of Analysis} (J. Koorevaar, ed.),
{\it AMS Proc.\ Symp.\ Pure Math.\ \bf XI} (1968), 360-369.
\medskip

[14] ---------, A Hilbert space of entire functions related to the
operational calculus, manuscript.
\medskip

[15] M.\ A.\ Rieffel, Non-commutative tori --- a case study of
non-commutative differentiable manifolds, {\it Contemp.\ Math.\ \bf 105}
(1990), 191-211.
\medskip

[16] K.\ Seip and R.\ Wallst\'en, Density theorems for sampling and
interpolation in the Bargmann-Fock space. II, {\it J.\ Reine Angew.\
Math.\ \bf 429} (1992), 107-113.
\medskip

[17] S.\ L.\ Woronowicz, Differential calculus on compact matrix
pseudogroups (quantum groups), {\it Comm.\ Math.\ Phys.\ \bf 122}
(1989), 125-170.
\bigskip

\noindent Department of Mathematics

\noindent Washington University

\noindent St.\ Louis, MO 63130

\noindent rr@math.wustl.edu

\noindent nweaver@math.wustl.edu
\end